\documentclass[12pt,reqno]{amsart}

\usepackage{latexsym}
\usepackage{amsfonts,amsmath,amssymb}
\usepackage{euscript}

\allowdisplaybreaks

\begin{document}
\title{An identity conjectured by Lacasse via the tree function}

\author{Helmut Prodinger}
\address{Department of Mathematics, University of Stellenbosch 7602,
Stellenbosch, South Africa}
\email{hproding@sun.ac.za}
\thanks{The  author was supported by an incentive grant of the NRF of South Africa.}

\begin{abstract}
A. Lacasse conjectured a combinatorial identity in his study of learning theory. Various people found independent proofs. Here is another one that is based on the study of the tree function, with links to Lamberts $W$-function and Ramanujan's $Q$-function. It is particularly short.
\end{abstract}

\maketitle

Let
\begin{align*}
\xi(n)&=\sum_{k=0}^n\binom nk \Bigl(\frac kn\Bigr)^k\Bigl(\frac {n-k}{n}\Bigr)^{n-k},\\
\xi_2(n)&=\sum_{k_1+k_2+k_3=n}\frac{n!}{k_1!k_2!k_3!}
 \Bigl(\frac {k_1}n\Bigr)^{k_1}\Bigl(\frac {k_2}n\Bigr)^{k_2}\Bigl(\frac {k_3}n\Bigr)^{k_3}.
\end{align*}
Then Lacasse~\cite{Lacasse} conjectured that $\xi_2(n)=\xi(n)+n$. There are now three independent
proofs~\cite{Younsi, CePeYa, Sun}. Here, we want to shed additional light on the matter, 
by using the \emph{tree function} (equivalent to Lambert's $W$-function~\cite{CoGoHaJeKn}) and linking the enumeration
to the celebrated $Q$-function of Ramanujan~\cite{FlGrKiPr}.

As Younsi has already pointed out, it is easier to work with
\begin{align*}
\alpha(n)&=n^n\xi(n)=\sum_{k=0}^n\binom nk k^k(n-k)^{n-k},\\
\beta(n)&=n^n\xi_2(n)=\sum_{k_1+k_2+k_3=n}\frac{n!}{k_1!k_2!k_3!}
 k_1^{k_1}k_2^{k_2}k_3^{k_3}.
\end{align*}

The tree function~\cite{GoJa} $y(z)$ is defined by $y=ze^y$ and possesses the expansion
\begin{equation*}
y(z)=\sum_{n\ge1}n^{n-1}\frac{z^n}{n!}.
\end{equation*}
Therefore
\begin{equation*}
\sum_{n\ge0}n^{n}\frac{z^n}{n!}=1+zy'(z)=1+\frac{y}{1-y}=\frac{1}{1-y}.
\end{equation*}
Consequently,
\begin{equation*}
\alpha(n)=n![z^n]\Bigl(\frac{1}{1-y}\Bigr)^2\quad\text{and}\quad
\beta(n)=n![z^n]\Bigl(\frac{1}{1-y}\Bigr)^3.
\end{equation*}
Now we compute the coefficients of a general power via Cauchy's integral formula:
\begin{align*}
[z^n]\Bigl(\frac{1}{1-y}\Bigr)^d
&=\frac1{2\pi i}\oint\frac{dz}{z^{n+1}}\Bigl(\frac{1}{1-y}\Bigr)^d\\*
&=\frac1{2\pi i}\oint\frac{dy(1-y)e^{-y}e^{y(n+1)}}{y^{n+1}}\Bigl(\frac{1}{1-y}\Bigr)^d\\
&=[y^n]\frac{e^{ny}}{(1-y)^{d-1}}\\
&=\sum_{k=0}^n\frac{n^{n-k}}{(n-k)!}\binom{k+d-1}{d-1}.
\end{align*}
Therefore
\begin{align*}
\alpha(n)=\sum_{k=0}^n\frac{n!}{k!}n^k,\qquad
\beta(n)=\sum_{k=0}^n\frac{n!}{k!}(n+1-k)n^k,
\end{align*}
and
\begin{align*}
\beta(n)-\alpha(n)&=\sum_{k=0}^n\frac{n!}{k!}(n-k)n^k\\
&=\sum_{k=0}^n\frac{n!}{k!}n^{k+1}-\sum_{k=1}^n\frac{n!}{(k-1)!}n^k=n^{n+1},
\end{align*}
as claimed.

Note that $\alpha(n)$ is linked to Ramanujan's celebrated $Q$-function~\cite{FlGrKiPr} via
\begin{equation*}
\alpha(n)=n^n\big(1+Q(n)\big).
\end{equation*}
\bibliographystyle{plain}

\end{document}